\documentclass[]{amsart}

\usepackage{geometry}
\usepackage{todonotes}
\usepackage{amsmath, amssymb, amsthm}
\usepackage[capitalize]{cleveref}

\numberwithin{equation}{section}

\newtheorem{theorem}{Theorem}[section]
\newtheorem{proposition}[theorem]{Proposition}

\newtheorem{lemma}[theorem]{Lemma}

\theoremstyle{definition}
\newtheorem{definition}[theorem]{Definition}
\newtheorem{remark}[theorem]{Remark}
\newtheorem{example}[theorem]{Example}

\DeclareMathOperator{\ad}{ad}
\DeclareMathOperator{\Rep}{\mathbf{Rep}}

\title{Semisimplifying Frank Lie algebras}
\author{Michiel Smet}

\begin{document}
	
	\maketitle
	
	\begin{abstract}
		The Frank Lie algebras are simple Lie algebras that only occur over fields of characteristic $3$. These come equipped with distinguished inner derivations that make them algebras in the category $\textbf{Rep}(\alpha_3)$. We apply the semisimplification functor to these Frank Lie algebras and obtain modular contact Lie superalgebras. 
		
		We also obtain a class of simple $J$-ternary algebras whose associated Jordan algebras are not simple.
	\end{abstract}

\section*{Introduction}

The recent construction of certain simple Lie superalgebras that only exist over fields of characteristic $3$ and $5$ developed by Kannan \cite{Kannan_2022} involving tensor categories, has led others to consider similar constructions for different kinds of algebras, e.g., \cite{Elduque2022, Elduque_Albert2025}. Most relevant for this article is a particular case of this construction that starts from a Lie algebra obtained from a $J$-ternary algebra over a field of characteristic $3$, see, e.g., \cite{Elduque_Cunha, Michiel2026}.
Namely, in these last two articles, it is observed that any Lie algebra $L$ coming from a $J$-ternary algebra $A$ is equipped with a canonical $\mathbb{Z}$-grading 
\[ L = L_{-2} \oplus L_{-1} \oplus L_0 \oplus L_1 \oplus L_2.\]
If $J$ is unital, there is a corresponding element $1 \in L_2$. With respect to the inner derivation $\ad 1$, the Lie algebra $L$ becomes a $\textbf{Rep}(\alpha_3)$-algebra, i.e., an algebra in the symmetric tensor category $\textbf{Rep}(\alpha_3)$. The semisimplification functor on $\textbf{Rep}(\alpha_3)$ maps $L$ to an operadic Lie superalgebra $\text{InDer}(A) \oplus A$.

In this article, we will look at the Frank Lie algebras $\mathcal{F}(n)$. These form a class of simple Lie algebras that only exist over fields of characteristic $3$ and can be constructed as subalgebras of the modular contact Lie algebras $K(3; (1,1,n))$. Frank \cite{Frank73} was the first to construct an algebra of this class. Using Strade's \cite{Strade} description of these algebras, we immediately obtain a $\mathbb{Z}$-grading and an inner derivation, analogous to the setting involving $J$-ternary algebras. We show that the image under the semisimplification functor is a modular contact Lie superalgebra $K(1,1; \underline{n})$.

In the final section, we also establish that the Frank Lie algebras can be obtained from certain $J$-ternary algebras.
Most notably, this construction yields a class of simple $J$-ternary algebras for which the associated Jordan algebra $J$ is not simple, contradicting \cite[Theorem 6.1]{Hein81}. We discuss the proof of this theorem and note that a simple $J$-ternary algebra either has a simple associated Jordan algebra $J$ or a \textit{degenerate} one in the sense of \cite{JordanAlgebren}.

\subsection*{Outline} Section 1 recalls the construction of the Lie (super)algebras in question. Section 2 introduces the semisimplification functor on $\textbf{Rep}(\alpha_3)$ and applies it to Frank algebras to obtain the modular contact Lie superalgebra. Section 3 introduces $J$-ternary algebras, defines the $J$-ternary algebra for a Frank algebra, and uses it to provide a counterexample to \cite[Theorem 6.1]{Hein81}.

\section{Basic definitions}	

Throughout, we work over a field $\Phi$ of characteristic $3$.

\subsection{Building blocks}
We use the algebras $\mathcal{O} = \mathcal{O}(1;\underline{n})$, $W = W(1;\underline{n})$ as defined by Strade \cite[2.1]{Strade}.
To be precise, $\mathcal{O}$ is the divided power algebra with basis elements $x^{(i)}$ for $0 \le i < 3^n$, and operation\footnote{This involves binomial coefficients, which are always integers and are thus defined over fields of arbitrary characteristic.}
\[ x^{(i)}x^{(j)} = \binom{i + j}{i} x^{(i+j)}\]
using $x^{(k)} = 0 $ for $k \ge 3^n$. 
This is a divided power algebra, since the elements $x^{(i)}$ mimic the behavior of $x^i/(i!)$.

The map $\partial(x^{(i)}) = x^{(i-1)}$, $i = 1,\ldots, 3^n - 1$ defines a derivation of $\mathcal{O}$.
We will often write $1$ for $x^{(0)}$ and $x$ for $x^{(1)}$.

Set $W$ to be the special derivations of $\mathcal{O}$, i.e., the set of all $f \partial $ with $f \in \mathcal{O}$. This forms a simple Lie algebra with operation
\[ [f\partial, g\partial] = (f\partial(g) - g\partial(f))\partial.\]
We can give $\mathcal{O}$ a second $W$-module structure using
\[ f\partial \cdot_{(\text{div})} g = f \partial(g) + \partial(f)g = \partial(fg).\]
We will write $\mathcal{O}_{(\text{div})}$ to indicate that we consider $\mathcal{O}$ with this second action.
For this second action, we have an irreducible module $\mathcal{O}'\subset \mathcal{O}_{(\text{div})}$, namely
\[ \mathcal{O}' = \bigoplus_{i = 0}^{3^n - 2} \Phi x^{(i)},\]
as proved by Strade \cite[Proposition 4.3.2]{Strade}.

\subsection{Modular contact Lie superalgebras}

The algebra $\mathcal{O}(1; \underline{n})$ has a super analogue  $\mathcal{O}(1,1; \underline{n}) \cong \mathcal{O}(1; \underline{n}) \otimes \Lambda(1)$ with $\Lambda(1)$ the Grassmann algebra on a single generator $\theta$.
Thus, in $\mathcal{O}(1,1; \underline{n})$ the elements are of the form $f + g \theta$, with $f$ the even part and $g \theta$ the odd part, for $f, g \in \mathcal{O}(1; \underline{n})$.
One defines the Witt superalgebra $W(1,1; \underline{n})$ as the Lie superalgebra formed by the elements
$ f \partial + g \partial_{\theta}$
with $f, g \in \mathcal{O}(1,1; \underline{n})$, where $\partial$ is the natural derivation on $\mathcal{O}(1; \underline{n})$, and $\partial_\theta$ represents the derivation associated to the generator of the Grassmann algebra, i.e., $\partial_{\theta} (f \theta^i) = i f \theta^{i - 1}$ for $i = 0, 1$.
The even part of $W(1,1; \underline{n})$ preserves the parity when it acts on $\mathcal{O}(1,1; \underline{n})$, while the odd part flips the parity.

In \cite{Liu05, Bai14}, a definition of the modular contact Lie superalgebra $K(1,1; \underline{n})$ is given. One way to define it is as the subalgebra of $W(1,1; \underline{n})$ generated by the image of the map
\[ D_K : \mathcal{O}(1,1; \underline{n}) \longrightarrow  W(1,1; \underline{n})\]
defined by
\[ D_K(x^{(i)}) = \theta x^{(i - 1)} \partial_\theta + 2x^{(i)} \partial  \]
and 
\[ D_K(x^{(i)} \theta ) = x^{(i)} \partial_\theta + x^{(i)} \theta  \partial. \]
In \cite{Bai14}, it is observed that $K(1,1; \underline{n}) \cong \mathcal{O}(1,1; \underline{n})$, where the latter Lie superalgebra comes equipped with the bracket
\[ [a,b] = D_K(a)(b) - 2 \partial(a) b.\]
This superalgebra is not simple, but the derived one is.
Now, we provide an even easier description of $K(1,1; \underline{n})$ and the derived algebra.

\begin{definition}
	We introduce the Lie superalgebra $W(1; \underline{n}) \oplus \mathcal{O}_{(\mathrm{div})}$. It has $W(1; \underline{n})$ as the even subalgebra, $\mathcal{O}_{(\mathrm{div})}$ as the odd part (on which the even part acts using the $(\mathrm{div})$-action), and bracket $[x,y] = - xy \partial $ for $x,y \in \mathcal{O}_{(\mathrm{div})}$. 
\end{definition}

\begin{lemma}
	\label{lem: contact}
	The Lie superalgebra $K(1,1; \underline{n})$ is isomorphic to $W(1; \underline{n}) \oplus \mathcal{O}_{(\mathrm{div})}$.
	Moreover, the derived algebra $K(1,1; \underline{n})^{(1)}$ is isomorphic to $W(1; \underline{n}) \oplus \mathcal{O}'$.
	\begin{proof}
		First, note that $[x^{(i)},x^{(j)}] = - x^{(i)} x^{(j - 1)} + x^{(j)} x^{(i - 1)}$, since the characteristic is $3$.
		Hence $x^{(i)} \mapsto - x^{(i)} \partial$ defines an isomorphism of the even parts of the Lie superalgebras.
		
		We also have $[x^{(i)}, x^{(j)} \theta] = - \theta( x^{(i - 1)} x^{(j)} + x^{(i)} x^{(j - 1)})$. So, $x^{(j)} \theta \mapsto x^{(j)}$ defines a morphism of the odd part of $K(1,1; \underline{n})$ onto $\mathcal{O}_{(\text{div})}$ compatible with the $W(1; \underline{n})$-action $f\partial \cdot_{(\text{div})} g = f \partial(g) + g \partial(f)$.
		
		In addition, one computes that $[x^{(i)}\theta, x^{(j)} \theta] = x^{(i)} x^{(j)}$. Hence, the bracket $\mathcal{O} \times \mathcal{O} \longrightarrow W(1; \underline{n})$ matches the one described in the lemma.

		Now, we determine the derived algebra.
		Note that $\mathcal{O}'$ is an irreducible submodule of $\mathcal{O}(1; \underline{n})$ and that $W(1; \underline{n})$ is simple. Hence, $K(1,1; \underline{n})^{(1)}$ certainly contains $W(1; \underline{n}) \oplus \mathcal{O}'$. Since $[f \partial, o] = \partial(fo)$ for $f, o \in \mathcal{O}$, it follows immediately that $[W,\mathcal{O}] \subset \mathcal{O}'$.
	\end{proof}
\end{lemma}

\begin{remark}
	Observe that $W(1; \underline{n}) \oplus \mathcal{O}'$ is a simple Lie superalgebra, given that the even part is a simple Lie algebra and the odd part is an irreducible module for that Lie algebra. Hence, $K(1,1; \underline{n})^{(1)}$ is simple.
\end{remark}

\subsection{Frank algebras}
\label{subsec Frank}

Consider a two dimensional vector space $U$ and a non-degenerate skew-symmetric form $\lambda$ on $U$.
Define $\phi : U \times U \longrightarrow \mathfrak{sl}(U)$ using
\[ \phi(u,u')(u'') = \lambda(u,u'') u' + \lambda(u',u'') u.\]

Strade \cite[4.4.D]{Strade} defines the Frank algebras $\mathcal{F}(n)$ as a certain class of Lie algebra with underlying vector space
\[ \mathcal{F}(n) = U \otimes \mathcal{O}_{(\text{div})} \oplus (\mathfrak{sl}(U) \otimes \mathcal{O}) \oplus \text{Id}_U \otimes W,\]
using the $2$-dimensional vector space $U$ and the algebras $\mathcal{O} = \mathcal{O}(1;\underline{n})$ and $W = W(1;\underline{n})$.

To describe the bracket, note that $\text{Id}_U \otimes W$ is a Lie subalgebra canonically isomorphic to $W$ and that $(\mathfrak{sl}(U) \otimes \mathcal{O})$ can be thought of as $\mathfrak{sl}(U)$ with scalars extended to $\mathcal{O}$. The algebra $\text{Id}_U \otimes W$ acts on $(\mathfrak{sl}(U) \otimes \mathcal{O})$, and both act naturally on $U \otimes \mathcal{O}_{(\text{div})}$.
For the bracket on $U \otimes \mathcal{O}_{(\text{div})} \times U \otimes \mathcal{O}_{(\text{div})}$, we use
\[ [ u \otimes f, v \otimes g] = \phi(u,v) \otimes (f \partial (g) - g \partial(f)) + \text{Id} \otimes \lambda(u,v) fg \partial.\]

We note that these algebras are simple. These algebras are related to the modular contact Lie algebras $K(3; (1, 1 , n))$, cfr. \cite[4.4.D]{Strade}.

\begin{remark}
	When Skryabin \cite{Skryabin92} gave a construction of this class of Lie algebras, the Lie superalgebra $K(1,1; \underline{n}) = W(1; \underline{n}) \oplus \mathcal{O}_{(\text{div})}$ did already appear in some sense. Namely, when trying to define the Frank algebra, Skryabin first verified the Jacobi identity (for the superalgebra) restricted to $\mathcal{O}_{(\text{div})}^3$ \cite[final equation page 393]{Skryabin92}, before proving the Jacobi identity on the Frank Lie algebra.
\end{remark}

\begin{remark}
	The construction yields a $\mathbb{Z}/2\mathbb{Z}$ grading with $(\mathfrak{sl}(U) \otimes \mathcal{O}) \oplus \text{Id}_U \otimes W$ as the subalgebra of even elements and $U \otimes \mathcal{O}_{(\text{div})}$ as the subspace of odd elements. Choosing a torus of $\mathfrak{sl}(U)$ and taking the corresponding weight spaces, refines this grading into a $\mathbb{Z}$-grading. On the even subspace, we have weight spaces with weight $0,$ $2$, and $-2$, and on the odd subspace we obtain weights $1$ and $-1$ (remark that these eigenspaces are really distinct over a field of characteristic $3$). 
\end{remark}

\section{Link between $K(1,1; \underline{n})$ and $\mathcal{F}(n)$}
In order to describe the connection with Frank algebras, we start by describing the semisimplification functor for the category $\Rep(\alpha_3)$.
For a more detailed description of a similar functor on $\Rep(\mathcal{C}_3)$, see \cite[section 3]{Ostrik2020} or for a very detailed description see \cite{Elduque2022}. For a lengthier description for $\Rep(\alpha_3)$, see \cite[section 3]{Kannan_2022}.

Recall that we work over a field $\Phi$ of characteristic $3$.

\subsection{Rep$(\alpha_3)$ and the semisimplification functor}

Set $\alpha_3(R) = \{ k \in R | k^3 = 0\}$ for commutative, associative, unital $\Phi$-algebras $R$. This forms an affine algebraic group with operation $(k,l) \mapsto k + l$. The coordinate algebra is $\Phi[t]/(t^3)$. The comultiplication and counit are determined by $\Delta(t) = t \otimes 1 + 1 \otimes t$ and $\epsilon(t) = 0$.

\begin{lemma}
	\label{lem: char rep algs}
	Representations of $\alpha_3$ on a vector space $V$ stand in bijection to maps $f : V \longrightarrow V$ such that $f^3 = 0$, using
	\[ k \cdot (1 \otimes_R v) = 1 \otimes_R v + k \otimes_R f(v) + k^2 \otimes_R f(f(v))/2\]
	for all $k \in \alpha_3(R)$.
	Moreover, if $V$ is an algebra, then it is an $\Rep(\alpha_3)$-algebra with respect to $f$ if and only if $f$ is a derivation.
	\begin{proof}
		An action of $\alpha_3$ on $V$ corresponds to a coaction $\delta$ of $\Phi[t]/(t^3)$ on $V$.
		We define $f,g,$ and $h$ using
		\[ \delta(v) = 1 \otimes h(v) + t \otimes f(v) + t^2 \otimes g(v).\]
		Using that $(\epsilon \otimes \text{Id}) \circ \delta = \text{Id}$ we conclude that $h(v) = v$.
		From $(\Delta \otimes \text{Id}) \circ \delta = (\text{Id} \otimes \delta) \circ \delta$ and comparison of scalars, we obtain
		\[ 2 t \otimes t \otimes g(v) = t \otimes t \otimes f(f(v))\]
		and 
		\[ 0 = t^2 \otimes t \otimes f(g(v)).\]
		This proves that each action of $\alpha_3$ is necessarily of the form described in the lemma. The converse, i.e., that each $f$ such that $f^3 = 0$ defines an action, is easily verified.

		For the second statement, we observe that
		$k \cdot (1 \otimes vw) = (k \cdot (1 \otimes v))(k \cdot (1 \otimes w)),$
		which is equivalent to $f(vw) = f(v)w + vf(w)$ and  $f(f(vw)) = f(f(v))w + vf(f(w)) + f(v)f(w)$. This proves that $V$ is a $\Rep(\alpha_3)$-algebra with respect to $f$ if and only if $f$ is a derivation.
	\end{proof}
\end{lemma}

The indecomposable objects of the category $\Rep(\alpha_3)$ are given by $(a)$ for $a = 1$, $2$, and $3,$ each with basis $\{e_i | 1 \le i \le a\}$ on which the $\Phi[f]/(f^3)$ representation is given by $f(e_i) = e_{i+1}$, writing $e_l = 0$ for $l > a$. The objects (2) and (3) are not simple, nor can be written as a direct sum of simple objects. Thus, the category is not semisimple.

We can semisimplify this category \cite[Proposition 2.13]{Ostrik2020} to obtain a semisimple category $\overline{\Rep(\alpha_3)}$.
This category is obtained by dividing out all negligible morphisms, i.e., morphisms $f:X \longrightarrow Y$ such that $\text{Tr}(fu) = 0$ for all $u : Y \longrightarrow X$. 
First, we remark that $\overline{(3)} \cong 0$ since $\text{Tr}(f) = 0$ for all $f : (3) \longrightarrow (3)$. Second, all morphisms $(2) \longrightarrow (1)$ and $(1) \longrightarrow (2)$ are negligible, so that $\text{Hom}(\overline{(i)},\overline{(j)}) = 0$ for $i \neq j$. Lastly, the identity maps $(i) \longrightarrow (i)$ are not negligible for $i = 1,2$, proving that $\text{End}(\overline{(i)}) \cong \Phi$.

\begin{lemma}
	The category $\overline{\Rep(\alpha_3)}$ is equivalent as a symmetric tensor category to the category of super vector spaces.
	\begin{proof}
		Both categories in question are generated by $2$ simple objects, say $O_0, O_1$, with $1$-dimensional endomorphism spaces as the only nontrivial homomorphism spaces, and all objects are isomorphic to direct sums of these $O_i$.
		Hence, it is sufficient to prove that the tensor products and braidings coincide, when expressed in terms of these $O_i$.
		
		For the supervectorspaces, we recall that this category is generated by a $1$-dimensional vector space $O_0 \cong L_0$ and a $1$-dimensional vector space $O_1 \cong L_1$ with $L_i \otimes L_j \cong L_{i +j}$ using mod $2$ indices,
		and braiding given by
		\[ L_i \otimes L_j \ni a \otimes b \longmapsto (-1)^{ij} b \otimes a \in L_j \otimes L_i.\]
		
		For $\overline{\Rep(\alpha_3)}$, we note that $\overline{(1)}$ can play the role of $O_0$ and $\overline{(2)}$ can play the role of $O_1$. The only nonobvious isomorphism is $O_1 \otimes O_1 \cong O_0$, which holds since
		\[ (2) \otimes (2)\cong (3) \oplus (1),\]
		in $\Rep(\alpha_3)$, using the decomposition of $(2) \otimes (2)$ into symmetric and antisymmetric tensors.
		
		We remark that $\Rep(\alpha_3)$ is trivially braided. This implies that only the braiding on $\overline{(2)} \otimes \overline{(2)}$ can be nontrivial. 
		Applying the $\text{Rep}(\alpha_3)$-braiding on the element $e = (v \otimes w - w \otimes v) \in (2) \otimes (2)$ yields $- e$, which proves that $f \circ c_{\overline{(2)} \otimes \overline{(2)}} = -f$ with $c_{A \otimes B} : A \otimes B \longrightarrow B \otimes A$ the braiding and $f$ the isomorphism $\overline{(2)} \otimes \overline{(2)} \longrightarrow \overline{(1)}$.
		This shows that we have the same braiding as for the category of super vector spaces.
	\end{proof}
\end{lemma}

\begin{remark}
	Since the semisimplification functor is a braided tensor functor, it preserves commutative diagrams that involve morphisms, the braiding, and tensor products.
	Consequently, a Lie $\Rep(\alpha_3)$-algebra maps to an (operadic) Lie superalgebra, i.e., an algebra $L_0 \oplus L_1$ that satisfies all linear axioms for Lie superalgebras but not necessarily the homogeneous Jacobi identity on $L_1$. 
\end{remark}

\subsection{Semisimplifying the Frank algebra} \label{sse:Frank}
Recall that the Frank algebra $\mathcal{F}(n)$ has as underlying vector space
\[ \mathcal{F}(n) = U \otimes \mathcal{O}_{(\text{div})} \oplus (\mathfrak{sl}(U) \otimes \mathcal{O}) \oplus \text{Id}_U \otimes W,\]
for a $2$-dimensional vector space $U$ and the algebras $\mathcal{O}, W$ as introduced before.
Recall that $U$ came equipped with a skew-symmetric form $\lambda$, from which we defined $\phi : U \times U \longrightarrow \mathfrak{sl}(U)$. Take $u ,v \in U$ such that $\lambda(u,v) = 1$. Note that $\phi(u,u)(v) = u$ and $\phi(u,u)(u) = 0$. Hence, we can think of $\phi(u,u)$ as the element
\[ e = \phi(u,u) = \begin{pmatrix}
	0 & 1 \\ 0 & 0
\end{pmatrix}\]
in $\mathfrak{sl}(U)$, using $u,v$ as the basis of $U$.

We look at the action of the algebra $\mathfrak{sl}(U)$ on the first tensor factor on $\mathcal{F}(n)$; we obtain a decomposition
\[ (2) \otimes \mathcal{O}_{(\text{div})} \oplus (3) \otimes \mathcal{O} \oplus (1) \otimes W\]
using the $(a)$ to denote standard irreducible $a$-dimensional $\mathfrak{sl}_2(U)$-module.
With respect to $d = \ad \left( e \otimes 1 \right)$, by using Lemma \ref{lem: char rep algs} and the fact that $(\ad e)^3 = 0$ on the involved $\mathfrak{sl}_2$-modules, we obtain a $\textbf{Rep}(\alpha_3)$-algebra. From the $\mathfrak{sl}(U)$-decomposition, we recover \[ (2) \otimes \mathcal{O}_{(\text{div})} \oplus (3) \otimes \mathcal{O} \oplus (1) \otimes W\]
as the decomposition into indecomposables $(1), (2),$ and $(3)$ of $\textbf{Rep}(\alpha_3)$.
So, applying the semisimplification functor yields a (operadic) Lie superalgebra with underlying module
\[ A'= W \oplus \mathcal{O}_{\text{(div)}},\]
as indicated in subsection 2.1. 

\begin{proposition}
	\label{prop}
	The semisimplification of $\mathcal{F}(n)$ with respect to $d$ is the contact superalgebra $K(1,1; \underline{n})$.
	\begin{proof}
		The preceding discussion shows that the semisimplification yields $W(1; \underline{n}) \oplus \mathcal{O}_{(\text{div})}$, equipped with a still to be determined product.
		We will show that the product matches the product of Lemma \ref{lem: contact}.
		
		The subspace $W$ is a subalgebra, since the semisimplification functor restricted to the subcategory of $\textbf{Rep}(\alpha_3)$ generated by (1) is an isomorphism to $\textbf{Vec}$.
		Any element $o \in \mathcal{O}_{(\text{div})}$, corresponds to 
		\[ \bar{o} : (2) \longrightarrow \mathcal{F}(n) : \alpha u + \beta v \mapsto (\alpha u + \beta v) \otimes o.\] 
		This shows that $W$ acts on $\mathcal{O}_{(\text{div})}$ using the (div)-action in the superalgebra, since $[w , \bar{o}(x)] = \overline{w(o)}(x)$ for $x \in (2)$.
		Finally, we use the isomorphism $\overline{(2)} \otimes \overline{(2)} \longrightarrow \overline{(1)}$ induced by the map $ (1) \oplus (3) \cong (2) \otimes (2)$ with $(1) \longrightarrow (2) \otimes (2) : \lambda \mapsto \lambda (u \otimes v - v \otimes u)$. Thus, for the (2)-indecomposables $\bar{o}$ and $\bar{o'}$, we obtain\footnote{Technically, one is interested in the result of the calculation below modulo $(2) \otimes \mathcal{O}_{(\text{div})} \oplus (3) \otimes \mathcal{O}$. However, since the result of the computation lies in $W \otimes (1)$, the space on which we should project, we can simply forget this. }
		\[ [\bar{o}, \bar{o'}] = [u \otimes o, v \otimes o'] - [v \otimes o, u \otimes o'] = 2 \lambda(u,v) o o' \partial = - o o' \partial. \qedhere\] 
	\end{proof}
\end{proposition}

\section{On J-ternary algebras}

Wolfgang Hein \cite{Hein81} defined a $J$-ternary algebra as a vector space $X$ equipped with a trilinear map $(x,y,z) \mapsto xyz$ such that
\begin{enumerate}
	\item $xy(uvw) - uv(xyw) = (xyu)vw + u(yxv)w,$
	\item $xyz - zyx = zxy - xzy$ 
\end{enumerate}
hold for all $x,y,z,u,v,w \in X$.

For such an algebra we can define 
\[ \langle x , y \rangle z = zyx - xyz.\]
Usually, $J$ denotes the span of all $\langle x,y \rangle$ for $x,y \in X$. We note that this is a Jordan algebra \cite[Theorem 1.1]{Hein81} with operation
\[ (a \bullet b) x = \frac{a(b(x)) + b(a(x))}{2}\]
for $a,b \in J$, $x \in X$.

\begin{example}
	\label{example Jtern}
	Consider $\mathcal{O}$ equipped with the trilinear map
	\[ (f,g,h) \mapsto L_{f,g} h = fg \partial h - hg \partial f.\]
	This defines a $J$-ternary algebra, as we will verify below. This product is closely linked to the Frank Lie algebras, which will become clearer when we verify the first axiom for J-ternary algebras.
	
	\begin{lemma}
		The space $\mathcal{O}$, equipped with the trilinear map $(f,g,h) \mapsto L_{f,g} h = fg \partial h - hg \partial f$, forms a $J$-ternary algebra.
		\begin{proof}
	We will use the product on the Frank Lie algebra, as introduced in subsection \ref{subsec Frank}, and take $u, v$ as in subsection \ref{sse:Frank} with $\lambda(u,v) = 1$.
	We first show that certain triple commutators can be understood in terms of the $J$-ternary algebra; afterwards we use the structure of the Frank Lie algebra to prove the first axiom for $J$-ternary algebras.
	Thus, in the Frank Lie algebra we first compute
	\begin{align*}
		 [[u \otimes f, v \otimes g], u \otimes h] & = [ \phi(u,v) \otimes (f \partial g - g \partial f) + \text{Id} \otimes fg \partial, u \otimes h] \\ & = - u \otimes h(f\partial g - g \partial f) + u \otimes (fg \partial \cdot_{\text{(div)}} h) \\
		 &= u \otimes (hg \partial f - hf \partial g) + u \otimes (fg \partial h + hf \partial g + hg \partial f)\\
		 &= u \otimes (gf \partial h - hg \partial f) = u \otimes L_{f,g} h.
	\end{align*}
	Similarly, one has 
		\begin{align*}
		[[v \otimes f, u \otimes g], v \otimes h] & = [ \phi(v,u) \otimes (f \partial g - g \partial f) - \text{Id} \otimes fg \partial, v \otimes h] \\ & = v \otimes h(f\partial g - g \partial f) - v \otimes (fg \partial \cdot_{\text{(div)}} h) \\
		&= - v \otimes (hg \partial f - hf \partial g) - v \otimes (fg \partial h + hf \partial g + hg \partial f)\\
		&= - v \otimes (gf \partial h - hg \partial f) = - v \otimes L_{f,g} h.
	\end{align*}
	
	From $[u \otimes f, v \otimes g]$ one can thus recover $L_{f,g}$ by acting on $u \otimes \mathcal{O}$, as indicated by the preceding computations. Thus, from
	\begin{align*} [[u \otimes x, v \otimes y], [u \otimes a, v \otimes b]] & = [[[u \otimes x, v \otimes y],u \otimes a],v \otimes b] - [u \otimes a,[[v \otimes y, u \otimes x], v \otimes b]] \\ & = [u \otimes L_{x,y} a, v \otimes b] + [u \otimes a, v \otimes L_{y,x} b] ,\end{align*}
	where we used the Jacobi identity and anti-commutativity in the first step and the preceding computations in the next step, we deduce
	\[ [L_{x,y}, L_{a,b}] = L_{(L_{x,y}a),b } + L_{a,(L_{y,x}b)}.\] Note that this is exactly the first condition on $J$-ternary algebras.
	The second condition is readily verified, using
	\[ L_{f,g}h - L_{h,g} f = 2(fg\partial h - hg \partial f)\]
	while
	\[ L_{h,f} g - L_{f,h} g = fh \partial g - gf \partial h - fh \partial g + gh \partial f = 2( gf \partial h - gh \partial f),\]
	using $-1 = 2$.	
\end{proof}
\end{lemma}
\end{example}

\begin{remark}
	One can always construct a Lie algebra from a $J$-ternary algebra. A general construction can be found in \cite{Hein81}. Applying this construction to the $J$-ternary algebra $\mathcal{O}$ as introduced in Example \ref{example Jtern}, yields a Frank algebra, as can be observed from the link between the commutators $[u \otimes f, v \otimes g]$ and the operators $L_{f,g}$ established in the previous lemma.
	
	Applying the construction of an operadic Lie superalgebra developed in \cite{Elduque_Cunha} or \cite{Michiel2026} to the $J$-ternary algebra (i.e., first construct the Lie algebra, then consider it as a $\textbf{Rep}(\alpha_3)$-algebra with respect to a distinguished derivation, and thereafter apply the semisimplification functor), yields the superalgebra $K(1,1; \underline{n})$ by Proposition \ref{prop}.
\end{remark}

\begin{lemma}
	\label{lem: J not simple}
	For the $J$-ternary algebra $\mathcal{O}$ as in Example \ref{example Jtern}, the associated Jordan algebra $J$ is isomorphic to the associative algebra $\mathcal{O}$ endowed with the ordinary Jordan product $x \bullet y = \frac{xy + yx}{2}$.
	Moreover, this algebra is not simple.
	\begin{proof}
		We already computed that
		\[ \langle f, h \rangle g =g (f\partial h - h \partial f).\]
		Thus $J$ can be identified with the span of certain operators $a \mapsto (f \partial h - h \partial f) a$ on $\mathcal{O}$. Observe the operators are (right or left) multiplications operators $R_b a = ab$ for $b = f \partial h - h \partial f$.
		This shows that $J$ can be identified with a subalgebra of the Jordan algebra $\mathcal{O}$, using $R_b \mapsto b$.
		To deduce that $J = \mathcal{O}$, recall first that $W(1; \underline{n})$ is simple so that its derived algebra coincides with $W(1; \underline{n})$. From 
		\[ J\partial = \langle f \partial h - h \partial f : f,h \in \mathcal{O}\rangle \partial = \langle [f\partial, h \partial] : f,h \in \mathcal{O} \rangle  = W(1,\underline{n})^{(1)} = W(1; \underline{n}) = \mathcal{O}\partial\]
		we conclude $J = \mathcal{O}$.
		
		The span of all $x^{(i)}$ with $i > 0$ forms an ideal of $\mathcal{O}$, thus $\mathcal{O}$ is not simple.
	\end{proof}
\end{lemma}

\begin{lemma}
	\label{lem: simple}
	The $J$-ternary algebra $\mathcal{O}$ is simple.
	\begin{proof}
	If $I$ is an ideal of the $J$-ternary algebra, the space $I\partial \le W$ would be an ideal of the simple algebra $W$, since 
	\[ \left(L_{f,1} h\right) \partial = [f\partial, h \partial]\]
	for all $h \in I$.
	Hence, $\mathcal{O}$ is a simple $J$-ternary algebra.
	\end{proof}
\end{lemma}

\begin{remark}
	\cite[Theorem 6.1]{Hein81} states that a $J$-ternary algebra $X$ is simple if and only if the associated Jordan algebra is simple and $\langle \cdot, \cdot \rangle$ is non-degenerate.
	Thus, Lemmas \ref{lem: J not simple} and \ref{lem: simple} contradict \cite[Theorem 6.1]{Hein81}. Indeed, $\mathcal{O}$ is a simple $J$-ternary algebra for which the associated Jordan algebra is not simple.
	
	Remark that the Jordan algebra $\mathcal{O}$ is degenerate in the sense of \cite[Chapter 1, Paragraph 8]{JordanAlgebren}. Whenever the associated Jordan algebra $J$ is non-degenerate, the (wrong part) of the proof of \cite[Theorem 6.1]{Hein81} can be recovered to prove that the Jordan algebra has to be simple. In the proof, Hein mistakenly claims that the nondegeneracy of $\langle \cdot, \cdot \rangle$ implies the non-degeneracy of the Jordan algebra $J$ by referring to section 4. In that section the non-degeneracy of $\langle \cdot, \cdot \rangle$ is used to establish the equivalence between the $J$-ternary algebra $X$ being non-degenerate and $J$ being non-degenerate; if the Killing form of the associated Lie algebra is non-degenerate, $X$ is also shown to be non-degenerate. Neither of these facts guarantees that $J$ is non-degenerate.  
	
	Given that plenty of simple Lie algebras have a non-degenerate Killing form, the $J$-ternary algebras $\mathcal{O}(1;\underline{n})$ might well form the sole class of simple degenerate $J$-ternary algebras over fields of characteristic different from $2$.
\end{remark}
\bibliographystyle{alpha}
\bibliography{bib}

\end{document}